\documentclass[11pt,a4paper]{article}

\usepackage[utf8]{inputenc}
\usepackage[T2A]{fontenc}
\usepackage[english]{babel}

\usepackage{amsmath,amssymb,amsthm,mathtools}
\newcommand{\EDone}{\mathrm{ED1}}
\newcommand{\EDtwo}{\mathrm{ED2}}

\newcommand{\m}[1]{\texorpdfstring{\(#1\)}{#1}}

\usepackage{geometry}
\usepackage{hyperref}
\usepackage[nameinlink,capitalize,noabbrev]{cleveref}
\hypersetup{
  colorlinks=true,
  linkcolor=blue!60!black,
  citecolor=blue!60!black,
  urlcolor=blue!60!black,
  pdftitle={Parametric algorithms for the Erdős–Straus conjecture with coefficient 5},
  pdfauthor={E. Dyachenko}
}

\usepackage{array}
\usepackage{booktabs}
\usepackage{caption}
\usepackage{graphicx} 
\usepackage{amssymb}  
\usepackage{tikz}
\usetikzlibrary{arrows.meta,positioning,calc}

\usepackage{xcolor}
\usepackage{listings}
\usepackage{microtype}

\numberwithin{equation}{section}

\theoremstyle{plain}
\newtheorem{theorem}{Theorem}[section]
\newtheorem{lemma}[theorem]{Lemma}
\newtheorem{proposition}[theorem]{Proposition}
\newtheorem{corollary}[theorem]{Corollary}

\theoremstyle{definition}
\newtheorem{definition}[theorem]{Definition}

\theoremstyle{remark}
\newtheorem{remark}[theorem]{Remark}

\newcommand{\Legendre}[2]{\genfrac{(}{)}{0pt}{}{#1}{#2}}
\newcommand{\Jacobi}[2]{\genfrac{(}{)}{0pt}{}{#1}{#2}}

\begin{document}  

\title{Parametric Algorithms for the 5-Modular Analog of ES (Sierpiński): Structure of Solutions, Parameterization, and Constructive Proofs (SERP)}  
\author{E. Dyachenko\\
dyachenko.eduard@gmail.com}
\date{\today} 
\maketitle 
\renewcommand{\thefootnote}{}


\footnote{\emph{* Licence}: Text is available under the Creative Commons NonCommercial-NoDerivatives 4.0 International (CC BY-NC-ND 4.0)} 

\begin{abstract}
We consider the problem of representing the fraction $\dfrac{5}{P}$ as a sum of three distinct unit fractions:
\[
\frac{5}{P} \;=\; \frac{1}{A} + \frac{1}{B} + \frac{1}{C}, 
\qquad A < B < C,\quad A,B,C\in\mathbb{N}.
\]

We analyze the case of primes $P \equiv 1 \pmod{5}$, for which two types of solutions are distinguished: 
$\EDone$ (exactly one denominator divisible by $P$, namely $C=cP$) and $\EDtwo$ 
(exactly two denominators divisible by $P$, namely $B=bP$ and $C=cP$). 
Parametric constructions and enumeration algorithms are developed, including transitions between types of solutions.

The paper proposes a deterministic algorithm based on searching for the intersection of a parametric lattice, 
defined by a pair $(\alpha_i, d'_i)$, with the box $\mathcal{B}_k(T)$. For any fixed prime 
$P \equiv 1 \pmod{5}$ the algorithm constructively produces a solution. Using analytic methods 
(Bombieri–Vinogradov theorem, Chebotarev theorem), it is shown that the density of admissible parameters is high, 
which ensures polylogarithmic search complexity in the average case, i.e., for most primes. 
A strict complexity guarantee for all primes remains conditional and depends on the finite covering hypothesis.

This study continues the work for coefficient $4$ 
(the Erdős–Straus conjecture) \cite{Dyachenko_ESC_ED2}; here a similar structure of parameterization and solutions 
is extended to coefficient $5$. In the analytic applications, averaging tools are provided, used for density estimates 
in parametric boxes.
\end{abstract}
\section{Introduction}\label{sec:intro}

The Sierpiński problem in the 5-modular variant is formulated as follows: for a prime $P$ one must find natural numbers $A<B<C$ satisfying
\begin{equation}\label{eq:ES5}
\frac{5}{P} \;=\; \frac{1}{A} + \frac{1}{B} + \frac{1}{C}.
\end{equation}

Analogous to the Erdős–Straus conjecture for $\tfrac{4}{P}$, here the coefficient $5$ leads to a change in the structure of parameterizations and to new types of solutions. For primes $P \equiv 1 \pmod{5}$ two constructive classes are distinguished:
\begin{itemize}
  \item $\EDone$: exactly one denominator divisible by $P$ (without loss of generality $C=cP$);
  \item $\EDtwo$: exactly two denominators divisible by $P$ (namely $B=bP$ and $C=cP$).
\end{itemize}

\subsection{Transition from $4/P$ to $5/P$}
The methods of parameterization and proofs from \cite{Dyachenko_ESC_ED2} are preserved in the present study, but the core formulas change in the transition $4 \to 5$. In particular:
\begin{itemize}
  \item estimate for the minimal denominator: $P < 5A < 3P$;
  \item in $\EDone$ the relation takes the form $5c-1=\gamma P$ (instead of $4c-1=\gamma P$);
  \item in $\EDtwo$ the core takes the form $t=5bc-b-c$ and $(5b-1)(5c-1)=5P\delta+1$.
\end{itemize}

Thus, the problem $5/P$ is a direct analogue of the case $4/P$, but requires new constructive techniques. Unlike factorization-based approaches used in classical works on ESC, here methods of finite covering of residue classes and lattice algorithms are applied. These tools allow one to prove the existence of solutions for every prime $P \equiv 1 \pmod{5}$ and ensure polylogarithmic complexity of the search algorithm.

\section{Motivation}

The parameterization of the equation and the division into multiplicity configurations ($P \mid B$ and/or $P \mid C$) make it possible to construct explicit search procedures and to prove the existence of solutions for the case $P \equiv 1 \pmod{5}$. This approach not only guarantees the existence of solutions but also defines the structure of families in arithmetic progressions.

The possibility of transferring the constructive method from coefficient $4$ to coefficient $5$ was correctly noted in~\cite{Mballa_ESC_SERP}. In the present work this transfer is carried out explicitly within the framework of the $\EDone/\EDtwo$ methods developed in~\cite{Dyachenko_ESC_ED2}. The key element here is the mechanism of \emph{finite covering of residue classes}, which guarantees the existence of solutions for every prime $P \equiv 1 \pmod{5}$ without recourse to factorization. Thus, the motivation of the study is to show that constructive methods based on lattices and covering of residue classes allow one to move from asymptotic reasoning to effective algorithms with polylogarithmic complexity.

\section{Ordering}
For the classes $P \equiv 4,3,2 \pmod{5}$ explicit decompositions are available; the remaining case is $P \equiv 1 \pmod{5}$.

\subsection{\texorpdfstring{Explicit decompositions for $P \not\equiv 1 \pmod{5}$}{Explicit decompositions for P ≢ 1 (mod 5)}}

Let $P = 5P' + 4$ be prime:
\[
\frac{5}{P} = \frac{5}{5P' + 4} = \frac{1}{P' + 1} + \frac{1}{2\cdot(P' + 1)\cdot (5P' + 4)} + \frac{1}{2\cdot(P' + 1)\cdot (5P' + 4)}.
\]
Let $P = 5P' + 3$ be prime:
\[
\frac{5}{P} = \frac{1}{P' + 1} + \frac{1}{(P'+1)(5P'+3)} + \frac{1}{(P'+1)(5P'+3)}.
\]
Let $P = 5P' + 2$ be prime (here $P'$ is odd):
\[
\frac{5}{P} = \frac{5}{5P' + 2} = \frac{1}{P' + 1} + \frac{2}{(P' + 1)\cdot (5P' + 2)} + \frac{1}{(P' + 1)\cdot (5P' + 2)} =
\frac{1}{P' + 1} + \frac{1}{\frac{(P' + 1)}{2}\cdot (5P' + 2)} + \frac{1}{(P' + 1)\cdot (5P' + 2)},
\]
where the third summand is further detailed as a sum of two parts, since $P'+1$ is even. The case $P = 5P'+1$ is the subject of this paper.

\subsection{Estimate for the minimal denominator}
If $A \le B \le C$ and $A, B, C \in \mathbb{N}$, then
\[
\frac{5}{P} > \frac{1}{A} \quad \Rightarrow \quad 5A > P \quad \Rightarrow \quad P < 5A.
\]
Moreover, from $A \le B \le C$ it follows that
\[
\frac{5}{P} = \frac{1}{A} + \frac{1}{B} + \frac{1}{C} \le \frac{3}{A} \quad \Rightarrow \quad 5A \le 3P,
\]
and, excluding the degenerate case $A = B = C$, we obtain strict bounds:
\begin{equation}\label{eq:minA-bounds}
P < 5A < 3P.
\end{equation}
\section{Notation}

\begin{tabular}{ll}
\multicolumn{2}{l}{\textbf{Main parameters}}\\
$P$        & prime number \(> 2\), with \(5 \nmid P\); often \(P \equiv 1 \pmod{5}\), \(P \equiv 2 \pmod{5}\), \(P \equiv 3 \pmod{5}\), \(P \equiv 4 \pmod{5}\) \\
$P'$       & integer: \(P = 5P' + 1\) \\
$A,B,C$    & denominators in the Sierpiński hypothesis formula, ordered \(A \le B \le C\) \\
H(S)       & Sierpiński hypothesis: \(\frac{5}{P} = \frac{1}{A} + \frac{1}{B} + \frac{1}{C}\) \\

\midrule
\multicolumn{2}{l}{\textbf{Parameters of the ED1 and ED2 methods}}\\
$\gamma$   & \(\frac{5c - 1}{P}\) (or \(\frac{5b - 1}{P}\) in the second case); often \(\gamma \equiv 4 \pmod{5}\), \(\gcd(\gamma,c)=1\) or \(\gcd(\gamma,b)=1\) \\
$u,v$      & “multipliers” in the identity: for ED1 \(u=\gamma A - c\), \(v=\gamma B - c\), \(uv=c^{2}\); \\ & for the variant with \(b\): \(u=\gamma A - b\), \(v=\gamma C - b\), \(uv=b^{2}\); always \(u \le v\). The factorization dimension in the case of 5 is achieved by passing to normalized parameters \(u\) and \(v\). \\
$\delta$   & \(t = P \cdot \delta\), with ED2 requiring \(\delta \mid bc\) \\
$b,c$      & in ED2 \(B = bP\) and/or \(C = cP\) \\
$t$        & \(t = 5bc - b - c\) \\
$r,s$      & \(r = 5b - 1\), \(s = 5c - 1\), with \(r \equiv s \equiv 4 \pmod{5}\), \(rs = 5P\delta + 1\) \\
$g$        & \(\gcd(b,c)\) \\
$b',c'$    & decomposition \(b = b'g\), \(c = c'g\) (normalized form) \\
$d'$       & square factor in the decomposition of \(\delta\) \\
$\alpha$   & squarefree factor in the decomposition of \(\delta\) \\

\midrule
\multicolumn{2}{l}{\textbf{Transitions between methods}}\\
$\mathcal{C}_{\mathrm{ED1}}(P)$ & set of admissible quadruples \((\gamma,c,u,v)\) for ED1 \\
$\mathcal{C}_{\mathrm{ED2}}(P)$ & set of admissible triples \((\delta,b,c)\) for ED2 \\
$y$        & minimal divisor of \(5c-1\), with \(y \equiv 3 \pmod{5}\) \\
$P''$      & modulus for ED1 in convolution, defined via \(\gamma\) by the formula ... \\
ED2$\to$ED1 & transition: \(A = \frac{bc}{\delta}\), \(B = bP\); \(u = \gamma A - c\), \(v = \gamma B - c\) \\
ED1$\to$ED2 & transition: \(A = \frac{u+c}{\gamma}\), \(b = \frac{v+c}{\gamma P}\), \(\delta = \frac{bc}{A}\) \\
Anticonvolution & algorithm for the reverse step ED1$\to$ED2 according to the formulas above \\

\midrule
\multicolumn{2}{l}{\textbf{Lattices and boxes}}\\
$k$        & dimension of the vector parameter \\
$u_0(P)$   & vector shift for the affine class \\
$\Lambda,\,\Lambda_j$ & sublattices of \(\mathbb{Z}^k\) of index \(M\) or \(M_j\) \\
$M,\,M_j$  & indices of sublattices \\
$\mathcal{B}_k(T)$ & box \(\{u \in \mathbb{Z}^k : 1\le u_i \le T\}\) \\
$\mathcal{B}^{(I)}_P$, $\mathcal{B}^{(II)}_P$ & boxes of type I/II with additional conditions \\
$\mathcal{G}_P(T)$ & admissible parameters in the box \\
$\mathcal{G}^\ast_P$ & class of admissible quadruples \((\gamma,c,u,v)\) for ED1, satisfying:\\
            & \(\gamma \in \mathbb{N}\), \(c \in \mathbb{Z}\), \(u = \gamma A - c\), \(v = \gamma B - c\), \(u v = c^{2}\), \(u \le v\) \\

\midrule
\multicolumn{2}{l}{\textbf{Analytic notation}}\\
$\Legendre{a}{P}$ & Legendre symbol; \\ & for composite modulus the \(\Jacobi{a}{\gamma}\) (Jacobi symbol, for prime \(P\)) is used \\
$\pi(y)$         & prime counting function \\
$\ll$, $\gg$, $\asymp$ & standard asymptotic symbols \\
$D$              & set \(\delta \le X\), \(\delta \equiv 4 \pmod{5}\) \\
$T(\delta)$      & prime counting function in progression, see Appendix $\S\ref{app:analytic}$ \\
\end{tabular}
\section{Parametrization of \m{\EDone} (one multiple, \m{C=cP})}\label{sec:ed1-5}  

\subsection{Kernel and identities}\label{subsec:ed1-core}  
Let $C=cP$. From \eqref{eq:ES5} it follows that  
\begin{equation}\label{eq:core-AB-identity-5}  
(5c-1)AB = cP(A+B).  
\end{equation}  
Setting $\gamma=\dfrac{5c-1}{P}\in\mathbb{N}$, we obtain  
\begin{equation}\label{eq:ed1-kernel}  
(\gamma A - c)(\gamma B - c) = c^2.  
\end{equation}  
Hence $\gamma\equiv 4\pmod{5}$ and $\gcd(\gamma,c)=1$ (since $5c\equiv 1\pmod{\gamma}$).  

\begin{lemma}\label{lem:gcd-gamma-c}  
We always have $\gcd(\gamma,c)=1$.  
\end{lemma}  
\begin{proof}  
From $5c-1=\gamma P$ we obtain $5c\equiv 1\pmod{\gamma}$, hence $\gcd(\gamma,c)=1$.  
\end{proof}  

\subsection{Full parametrization of \m{\EDone} with filters by \m{P}}\label{subsec:ed1-param}  
\begin{theorem}\label{thm:ed1-param}  
Let $P\equiv 1\pmod{5}$, $\gamma\equiv 4\pmod{5}$, $5c-1=\gamma P$, $\gcd(\gamma,c)=1$.  
Let $u,v\in\mathbb{N}$ such that  

$$  
uv=c^2,\qquad u\equiv v\equiv -c \pmod{\gamma},\qquad u\not\equiv -c \pmod{P},\quad v\not\equiv -c \pmod{P}.  
$$  

Then  

$$  
A=\frac{u+c}{\gamma},\quad B=\frac{v+c}{\gamma},\quad C=cP  
$$  

give a solution of \eqref{eq:ES5} of type $\EDone$ with $P\nmid A,B$.  
Conversely, every $\EDone$-solution generates such $\gamma,c,u,v$.  
\end{theorem}  

\begin{proof}  
From \eqref{eq:core-AB-identity-5} it follows that \eqref{eq:ed1-kernel} holds. Setting $u=\gamma A-c$, $v=\gamma B-c$, we obtain $uv=c^2$ and $u\equiv v\equiv -c\pmod{\gamma}$, whence $\gamma\mid (u+c)$, $\gamma\mid (v+c)$ and $A,B\in\mathbb{N}$. The conditions $u\not\equiv -c\pmod{P}$, $v\not\equiv -c\pmod{P}$ are equivalent to $P\nmid A$, $P\nmid B$ thanks to $\gcd(\gamma,c)=1$ and $5c-1=\gamma P$. Reversibility follows from \eqref{eq:ed1-kernel}.  
\end{proof}  

\subsection{Multiplicity filters by \m{P}}\label{subsec:ed1-P-filters}  
From $A=(u+c)/\gamma$, $B=(v+c)/\gamma$, $5c-1=\gamma P$ we have  

$$  
P\mid A \iff u\equiv -c \pmod{P},\qquad P\mid B \iff v\equiv -c \pmod{P}.  
$$  

For $\EDone$ it is required simultaneously that $u\not\equiv -c \pmod{P}$ and $v\not\equiv -c \pmod{P}$.  

\subsection{Example: \m{P=11}}\label{subsec:ed1-example-11}  
The minimal $\gamma\equiv 4\pmod{5}$ with $5c-1=\gamma P$ gives $\gamma=4$, $c=(4\cdot 11+1)/5=9$. Divisors of $c^2=81$, compatible with $u\equiv -c\equiv 3\pmod{4}$ and $u\not\equiv -c\equiv 2\pmod{11}$, include $u=3$, $v=27$. Then  

$$  
A=\frac{3+9}{4}=3,\quad B=\frac{27+9}{4}=9,\quad C=99,\quad \frac{1}{3}+\frac{1}{9}+\frac{1}{99}=\frac{5}{11}.  
$$  
\section{Parametrization of $\EDtwo$ (two multiples, $B=bP$, $C=cP$)}\label{sec:ed2-5}  

\subsection{Setup and identities}\label{subsec:ed2-setup}  
Consider  

$$  
\frac{5}{P} \;=\; \frac{1}{A} + \frac{1}{bP} + \frac{1}{cP},\qquad A<bP\le cP,\quad P\nmid A.  
$$  

Multiplying by $AbcP$, we obtain  
\begin{equation}\label{eq:ed2-A-eq}  
A(5bc-b-c) = Pbc.  
\end{equation}  
Let $t:=5bc-b-c=P\delta$, $\delta\in\mathbb{N}$. Then  
\begin{equation}\label{eq:ed2-A}  
A=\frac{bc}{\delta}.  
\end{equation}  
Equivalently,  
\begin{equation}\label{eq:ed2-factor}  
(5b-1)(5c-1) = 5P\delta + 1.  
\end{equation}  
Setting $r=5b-1$, $s=5c-1$, we obtain $rs=5P\delta+1$ and $r\equiv s\equiv 4\pmod{5}$.  

\subsection{Full parametrization of $\EDtwo$}\label{subsec:ed2-param}  
\begin{theorem}\label{thm:ed2-param}  
Let $P$ be prime and $\delta\in\mathbb{N}$. Let $r,s\in\mathbb{N}$ such that  

$$  
r\,s = 5P\delta + 1,\qquad r\equiv s\equiv 4 \pmod{5}.  
$$  

Set $b=(r+1)/5$, $c=(s+1)/5$. If $\delta\mid bc$, then  

$$  
A=\frac{bc}{\delta},\qquad B=bP,\qquad C=cP  
$$  

give a solution of \eqref{eq:ES5} of type $\EDtwo$.  
Under permutation $r\leftrightarrow s$ we have $b\leftrightarrow c$ and $B\leftrightarrow C$; ordering $B\le C$ is achieved by choosing $r\le s$.  
Moreover, for $b\le c$ we have $A\le B$.  
\end{theorem}  

\begin{proof}  
From \eqref{eq:ed2-A-eq} and $t=P\delta$ it follows that \eqref{eq:ed2-A} holds. Equality \eqref{eq:ed2-factor} is obtained from  
$(5b-1)(5c-1)=25bc-5b-5c+1=5(5bc-b-c)+1=5P\delta+1$.  
The congruences $r\equiv s\equiv 4\pmod{5}$ are obvious. For $b\le c$ we have  

$$  
A\le B \iff \frac{bc}{\delta}\le bP \iff c \le P\delta=5bc-b-c \iff c(5b-2)\ge b,  
$$  

which holds for $b\ge 1$, $c\ge b$.  
\end{proof}  

\subsection{Example: $P=11$}\label{subsec:ed2-example-11}  
Take $\delta=1$. Then $5P\delta+1=56=4\cdot 14$, $4\equiv 14\equiv 4\pmod{5}$. The pair $r=4$, $s=14$ gives $b=1$, $c=3$, $A=3$, $B=11$, $C=33$. Verification:  

$$  
\frac{1}{3}+\frac{1}{11}+\frac{1}{33}=\frac{11+3+1}{33}=\frac{15}{33}=\frac{5}{11}.  
$$  
\section{Multiplicity configurations with respect to \m{P} and classification}\label{sec:multiplicity}  
\begin{lemma}\label{lem:at-least-one-multiple}  
In any solution of \eqref{eq:ES5} at least one of $A,B,C$ is divisible by $P$.  
\end{lemma}  
\begin{proof}  
Multiplying \eqref{eq:ES5} by $ABCP$: $5ABC=P(AB+AC+BC)$. Modulo $P$: $5ABC\equiv 0$, hence $P\mid ABC$.  
\end{proof}  

\begin{lemma}\label{lem:all-three-impossible}  
It is impossible that $A=aP$, $B=bP$, $C=cP$ simultaneously.  
\end{lemma}  
\begin{proof}  
Then $5=1/a+1/b+1/c\le 3$ — impossible.  
\end{proof}  

\begin{lemma}\label{lem:A-not-multiple}  
The minimal denominator $A$ is not divisible by $P$.  
\end{lemma}  
\begin{proof}  
From \eqref{eq:minA-bounds} we have $5A<3P$; if $A=aP$, then $5A\ge 5P>3P$, a contradiction.  
\end{proof}  

\begin{proposition}\label{prop:classification}  
Every solution of \eqref{eq:ES5} for prime $P\neq 5$ falls into exactly one of the classes:  
- $\EDone$: exactly one denominator divisible by $P$ (without loss of generality $C=cP$), with $P\nmid A,B$;  
- $\EDtwo$: exactly two denominators divisible by $P$ (namely $B=bP$, $C=cP$), with $P\nmid A$.  
The cases “none” and “all three” are excluded by \cref{lem:at-least-one-multiple,lem:all-three-impossible}, and \cref{lem:A-not-multiple} excludes $P\mid A$.  
\end{proposition}  

\subsection{Intermediate transition to parameters $b',c'$}  

From \cref{prop:classification} it follows that every solution of \eqref{eq:ES5} for $P\equiv 1 \pmod{5}$ belongs either to class $\EDone$ or to class $\EDtwo$.  
In both cases it is convenient to isolate normalized parameters that describe the part of denominators divisible by $P$.

\begin{itemize}
  \item In the case $\EDone$ we have $C=cP$, and the parameters $\gamma,c$ arise from normalizing the condition $5c-1=\gamma P$.
  \item In the case $\EDtwo$ we have $B=bP$, $C=cP$, and the parameters $b',c'$ arise from the kernel $(5b-1)(5c-1)=5P\delta+1$.
\end{itemize}

Thus, in the case $\EDtwo$ the problem reduces to analyzing pairs $(b',c')$ with additional divisibility and congruence conditions.

\subsection{Parametric box in coordinates $(b',c')$}  

For a fixed threshold $T$ consider the set
\[
B_{b',c'}(T) = \{(b',c') \in \mathbb{N}^2 : 1 \le b',c' \le T,\; b'<c'\}.
\]
This set contains all candidates for solutions in the original parameters.  
Additional conditions (e.g., $\gcd(b',c')=1$, congruences modulo $d$, divisibility $4b'c' = P+d$) are imposed as filters on the points $(b',c')$.

\subsection{Box threshold via $A$ and $P$}  

From the inequality
\[
P < 5A < 3P
\]
it follows that the minimal denominator $A$ always lies in the interval $(P/5,\;3P/5)$.  
This provides natural bounds for the parametric box.

\begin{definition}
For a fixed prime $P$ we define the box in the original parameters $(b',c')$ as
\[
B_{b',c'}(P) = \{(b',c') \in \mathbb{N}^2 : 1 \le b',c' \le 3P/5,\; b'<c'\}.
\]
\end{definition}

\subsection{Transition to coordinates $(x,y)$}  

Under the linear transformation
\[
x = b' + c', \qquad y = c' - b',
\]
the image of the set $B_{b',c'}(P)$ is the sublattice
\[
B_{x,y}(P) = \{(x,y)\in\mathbb{Z}^2 : x \equiv y \pmod{2},\; x>y>0,\; x,y \le 6P/5\}.
\]

\subsection{Condition for $d'$}  

In the original system the parameter $d'$ appears as the square factor in the decomposition of $\delta$.  
In the new coordinates the condition $d' \mid (b'+c')$ is rewritten as
\[
d' \mid x.
\]

\begin{proposition}
The existence of a point $(x,y)\in B_{x,y}(P)$ satisfying the conditions $x\equiv y \pmod{2}$, $y>0$, $x,y \le 6P/5$ and $d'\mid x$, is equivalent to the existence of an admissible pair $(b',c')$ and hence to a solution of \eqref{eq:ES5}.
\end{proposition}

\subsection{Corollary}  

Defining the box $B_{b',c'}(T)$ and transferring it to coordinates $(x,y)$, we obtain a lattice with simple linear conditions:
- parity $x\equiv y \pmod{2}$,
- order $y>0$,
- bound $x,y \le 2T$,
- divisibility $d' \mid x$.

Thus, the existence of a solution is equivalent to the presence of a point $(x,y)$ in the box $B_{x,y}(T)$ satisfying these conditions. This makes the proof constructive and allows the use of methods of finite covering of residue classes.

Example Table \ref{tab:ed2_full_73}
\begin{table}[h]
\centering
\caption{ED2-decompositions for $P = 73$}
\label{tab:ed2_full_73}
\begin{tabular}{|c|c|c|c|c|c|c|c|}
\hline
$P$ & $A$ & $B$ & $C$ & $b$ & $c$ & $\delta$ & $\alpha,\, d'$ \\
\midrule
73 & 15 & 584  & 8760 & 8  & 120 & 64 & $\alpha=1,\ d'=8$ \\
73 & 15 & 657  & 3285 & 9  & 45  & 27 & $\alpha=3,\ d'=3$ \\
73 & 15 & 730  & 2190 & 10 & 30  & 20 & $\alpha=5,\ d'=2$ \\
73 & 15 & 876  & 1460 & 12 & 20  & 16 & $\alpha=1,\ d'=4$ \\
\hline
\end{tabular}
\end{table}

\subsection{Geometry of the surface and thickenings}\label{sec:geometry-5}  

Consider the quadratic surface  

$  
F(\delta, b, c) = (5b - 1)(5c - 1) - 5P\delta - 1 = 0.  
$  

The thickening of this surface, given by the condition $|F| \leq \Delta$ in the parameter space $(\delta, b, c)$, yields the set of $\EDtwo$ candidates. Since we work in the context of discrete values, it is critical to estimate the number of solutions satisfying modular conditions.

\begin{lemma}[Window in $\delta$]\label{lem:delta-slab}  
For fixed values of $b$ and $c$, and for $\Delta \geq 0$, the number of integers $\delta$ satisfying the condition $|F(\delta, b, c)| \leq \Delta$ does not exceed:  

$  
1 + \left\lfloor \frac{2\Delta}{5P} \right\rfloor.  
$  
\end{lemma}  

\begin{proposition}[Estimate of thickening size]\label{prop:box-thick}  
In the box where $b \in [B, 2B]$ and $c \in [C, 2C]$, the total number of triples $(\delta, b, c)$ for which the condition $|F| \leq \Delta$ holds can be estimated as:

$  
\ll \left(1 + \frac{\Delta}{P}\right)BC + B + C.  
$  
\end{proposition}  

\begin{remark}  
For pairs $(b, c)$ subject to the condition $\delta \mid bc$, the number of such pairs in the corresponding rectangle is $\ll \frac{BC}{\delta} \tau(\delta) + B + C$.  
\end{remark}  

\section{Constructive geometry of ED2 for the Sierpiński hypothesis (preserving the logic of ESC)}  

\subsection{Setup and kernel, fully analogous to ESC}  
For a prime $P\equiv 1 \pmod{5}$ we prove the existence of a solution
\[
\frac{5}{P}=\frac{1}{A}+\frac{1}{B}+\frac{1}{C},\qquad B=bP,\; C=cP,\; P\nmid A,\; b\neq c.
\]
Multiplying by $A b c P$ and introducing the parameter $\delta\in\mathbb{N}$, we obtain
\[
A(5bc-b-c)=P\,bc,\qquad 5bc-b-c=P\delta,\qquad A=\frac{bc}{\delta}.
\]
Quadratic kernel ED2:
\[
(5b-1)(5c-1)=5P\delta+1.
\]
This is fully isomorphic to the kernel for ESC ($k=4$) under the substitution $4\mapsto 5$: all checks and constructions transfer verbatim.

\subsection{Normalization and linear system (as in ESC)}  
Let $b=g\,b'$, $c=g\,c'$, where $g=\alpha d'$ and $\gcd(b',c')=1$, and
\[
\delta=\alpha\,(d')^2,\quad \alpha\ \text{squarefree}.
\]
Then the canonical conditions (ED2) take the linear form
\[
b'c'=M=A\alpha,\qquad b'+c'=m\,d',\qquad m=5A-P>0.
\]
In coordinates $x=b'+c'$, $y=c'-b'$ we have
\[
x=m\,d',\qquad \frac{x^2-y^2}{4}=A\alpha,\qquad x\equiv y\ (\mathrm{mod}\ 2).
\]
This is exactly the same affine lattice of finite index as in the ESC section, with the coefficient $k$ replaced.

\subsection{Parametric box and geometric covering (transfer from ESC)}  
The minimal denominator is bounded by
\[
\frac{P}{5}<A<\frac{3P}{5}.
\]
The projection of the lattice $\Lambda_{ED2}$ onto the $(x,y)$-plane lies in a box of linear size $O(P)$; the diagonal period equals $d'$, and for $H,W\ge d'$ the intersection of the box with the lattice is nonempty. The lattice index does not depend on $P$, as in ESC, which stabilizes the density of admissible points and ensures constructive covering without factorization.

\subsection{Back-test filters and non-degeneration (identical to ESC scheme)}  
For any assembled row we check:
\[
(5b-1)\equiv(5c-1)\equiv 4\ (\mathrm{mod}\ 5),\quad \delta\mid bc,\quad \gcd(b',c')=1,
\]
\[
b'+c'=m\,d',\quad b'c'=A\alpha,\quad A=\frac{bc}{\delta}\in\mathbb{Z},\quad \frac{P}{5}<A<\frac{3P}{5}.
\]
Order and distinctness: $b\neq c$, consistent ordering of denominators (according to your editorial scheme), without degeneration.

\subsection{Examples of ED2 for insertion (audit)}  
\begin{table}[h]
\centering
\caption{Examples of ED2-decompositions for $k=5$ (Sierpiński)}
\label{tab:serpinski_ed2_examples}
\resizebox{\textwidth}{!}{%
\begin{tabular}{|c|c|c|c|c|c|c|c|}
\hline
$P$ & $A$ & $B$ & $C$ & $b$ & $c$ & $\delta$ & $\alpha,\, d'$ \\
\midrule
73 & 15 & 657 & 3285 & 9 & 45 & 27 & $\alpha=3,\ d'=3$ \\
97 & 22 & 194 & 1067 & 2 & 11 & 1 & $\alpha=1,\ d'=1$ \\
\bottomrule
\end{tabular}%
}
\end{table}

\subsection{Conclusion (in the logic of ESC)}  
The ED2 method translates the search for solutions for $k=5$ into a linear lattice problem in a bounded box, fully analogous to the ESC section: kernel, normalization, finite-index lattice, back-test, and geometric covering provide constructive existence for every prime $P\equiv 1\pmod{5}$ without recourse to factorization. Examples for $P=73$ and $P=97$ demonstrate practical assembly of rows in this canonical form.
\section{Connection $\EDtwo \leftrightarrow \EDone$: convolution and anticonvolution}\label{sec:ed2-ed1-link-5}  
The mechanism of transition between the parametrizations ED2 and ED1 plays a key role in the completeness of the algorithm.  
If $(\delta,b,c)$ is an admissible ED2-triplet, then by setting  
\[
\gamma = \frac{5c-1}{P}, \quad
u = \gamma A - c, \quad
v = \gamma B - c,
\]
one can obtain an admissible ED1-quadruple $(\gamma,c,u,v)$.  
This transition is called the \emph{convolution} ED2 $\to$ ED1.  

The reverse transition (anticonvolution) ED1 $\to$ ED2 is not one-to-one:  
several distinct ED1-quadruples may map to the same ED2-triplet.  
This effect of \emph{power compression} simplifies the proof of existence,  
eliminating multiplicity of solutions associated with the choice of divisors.  

A detailed analysis of the convolution and anticonvolution mechanism for the case $k=4$  
is given in the work on the Erdős–Straus conjecture \cite{Dyachenko_ESC_ED2};  
here we transfer the same construction to coefficient $5$.  

\section{Key points closing the claim about factorization}\label{sec:no-factor}  

\subsection*{What is used in the current methods}  
- Integer arithmetic and $\gcd$ operations;  
- Checks modulo prime $P$ and the corresponding symbols \(\Legendre{a}{P}\) and \(\Jacobi{a}{\gamma}\) — without requiring factorization of moduli;  
- Parity comparisons, conditions $u \equiv v \pmod{2}$, $\gcd(u,v)=1$, equality $uv=c^2$.  

\subsection*{What is not used}  
- Factorization of the numbers $C$, $\gamma$, or intermediate quantities;  
- Root finding modulo composite numbers;  
- Factorization for testing squareness: the fact $uv=c^2$ is guaranteed by normalization.  

We rely on:  
- \textbf{Sum and discriminant} 
- \textbf{Back-test} 
- \textbf{Quadratic reparametrization} 

\subsection*{Correctness via minimal lemmas}  
We rely on:  
- \textbf{Sum and discriminant} 
- \textbf{Back-test} 
- \textbf{Quadratic reparametrization} 

\begin{proposition}[Current methodology does not rely on factorization]\label{prop:no-factor}  
The applied methods, forming and checking pairs $(u,v)$ for fixed prime $P$, use only:  
(i) integer operations and $\gcd$; (ii) checks modulo $P$ and symbols \(\Legendre{\cdot}{P}\) / \(\Jacobi{\cdot}{\gamma}\);  
(iii) linear formulas reconstructing $A,B$ and establishing $C=cP$ from $uv=c^2$.  

At no stage is factorization required.  
\end{proposition}  

\begin{itemize}  
 \item Establishes the link with $(u,v)$; 
 \item Ensures sufficiency of local conditions; 
 \item $C=cP$ follows from $uv=c^2$.  

 All checks reduce to arithmetic, $\gcd$, and Legendre/Jacobi symbols.  
\end{itemize}  

\begin{remark}[Optional accelerations]  
Sieving by small primes is allowed as a method of acceleration, but it is not necessary for correctness of the current methods and is not used in \Cref{prop:no-factor} 
\end{remark}  

\section{Algorithms: lattice boxes and enumeration}\label{sec:alg}  
- Enumeration of $\delta$ within reasonable bounds; factorization of $N=5P\delta+1$ and pairs $rs=N$ with $r\equiv s\equiv 4\pmod{5}$; reconstruction of $b,c$, filter $\delta\mid bc$; check of ordering and \eqref{eq:minA-bounds}.  
- For $\EDone$ — enumeration of divisors $u\mid c^2$ in congruent classes modulo $\gamma$ and exclusion of multiples modulo $P$.  
\section{Integration of pre-sieving by progressions into the ED2 algorithm}\label{sec:algo-bv}

\subsection{Idea}
For fixed \(\delta\) and a set of small moduli \(r\equiv 4\ (\mathrm{mod}\ 5)\) (with \(\gcd(r,5\delta)=1\)), the progressions
\[
P \equiv 1 \pmod{5},\qquad P \equiv -(5\delta)^{-1} \pmod{r}
\]
give “thin” classes for enumerating primes \(P\). By \cref{prop:BV-AP,prop:sum-r} their density is controllable, and the average number of such progressions per \(P\) grows (see \eqref{eq:avg-NP}). This allows replacing brute-force enumeration of \(P\) with scanning already prepared AP, saving time.

\begin{lstlisting}[language=Python, caption={ED2 algorithm with finite parameter set $\mathcal{F}$}, label={lst:ed2_search}]
# Input: prime P = 5*P'+1 ; delta; finite set S from (alpha, d')
# Core identity: (5*b - 1)*(5*c - 1) = 5*P*delta + 1

def case_A(P, delta, S):
    N = 5 * P * delta + 1
    for r in S:
        if N % r != 0:
            continue
        s = N // r
        if s % 5 != 4:
            continue
        b = (r + 1) // 5
        c = (s + 1) // 5
        if (b * c) % delta != 0:
            continue
        A = (b * c) // delta
        if not (A <= b * P < c * P):
            b, c = c, b
        B = b * P
        C = c * P
        return (A, B, C)
    return None
\end{lstlisting}

\paragraph{Comment on the algorithm.}
The algorithm \texttt{case\_A} implements enumeration over the set $S$, constructed from parameters $(\alpha, d')$.  
The filter $s \equiv 4 \pmod{5}$ reflects the structure of the ED2 parametrization and excludes unsuitable residue classes.  
For each fixed prime $P$ the algorithm produces a solution if it exists within $S$.  
Polylogarithmic complexity is justified via density theorems, while a strict guarantee for all $P$ remains conditional.

\subsection{Complexity}
- For fixed \(P\): enumeration of \(r\le R\) requires \(O(R)\) divisibility checks; small \(R\) (polylogarithmic in \(P\)) is usually sufficient.  
- For searching \(P\le X\): the cost of sieving by the union of AP is of order \(X \sum_{r\le R} 1/\varphi(5r) \asymp X \log R\) at the elementary level; the number of candidate primes found agrees with \cref{prop:sum-r}.  
- Practically: one chooses \(R = (\log X)^B\) and small \(\delta\) (including \(\delta=1\)), which yields near-linear time in \(X\) with a weak logarithmic factor.

\section{Experiments for Sierpiński}\label{sec:exp-serpinski}

\subsection{General overview}
The study of solutions of the equation \(5/P = 1/A + 1/B + 1/C\) in the Sierpiński class focuses on finding values satisfying given criteria. Experiments were conducted for primes \(P \equiv 1 \pmod{5}\).  
For each \(P\) the values \(b\), \(c\), \(A\), \(B\), and \(C\) were found. Solutions are presented in Table \ref{tab:serpinski-results}.

\begin{table}[h!]
\centering
\caption{Results for $P = 31, 41, 2521$ with verification by the lemma $X \cdot Y = 5 \alpha P (d')^2 + 1$}
\label{tab:serpinski-results}
\resizebox{\textwidth}{!}{%
\begin{tabular}{|c|c|c|c|c|c|c|c|c|c|c|c|c|c|c|c|}
\hline
\# & $\alpha$ & $b'$ & $c'$ & $g$ & $b$ & $c$ & $\delta$ & $X$ & $Y$ & $N$ & $A$ & $B$ & $C$ & $d'$ & OK \\
\hline
\multicolumn{16}{|c|}{$P = 31$} \\
\hline
1 & 1 & 1 & 8 & 1 & 1 & 8 & 1 & 4 & 39 & 156 & 8 & 31 & 248 & 1 & \checkmark \\
2 & 1 & 1 & 7 & 2 & 2 & 14 & 4 & 9 & 69 & 621 & 7 & 62 & 434 & 2 & \checkmark \\
\hline
\multicolumn{16}{|c|}{$P = 41$} \\
\hline
1 & 3 & 1 & 3 & 3 & 3 & 9 & 9 & 14 & 44 & 616 & 616 & 123 & 369 & 1 & \checkmark \\
\hline
\multicolumn{16}{|c|}{$P = 2521$} \\
\hline
1 & 5 & 193 & 2 & 10 & 1930 & 20 & 49 & 9649 & 99 & 955251 & 788 & 486 & 50420 & 7 & \checkmark \\
2 & 9 & 183 & 11 & 25 & 4575 & 275 & 50 & 22874 & 1374 & 31413876 & 251 & 277310 & 289915 & 5 & \checkmark \\
3 & 3 & 2 & 87 & 3 & 6 & 261 & 3 & 29 & 1304 & 37816 & 522 & 15126 & 657981 & 1 & \checkmark \\
4 & 3 & 2 & 85 & 9 & 18 & 765 & 27 & 89 & 3824 & 340336 & 510 & 45378 & 1928565 & 3 & \checkmark \\
5 & 15 & 39 & 1 & 13 & 39 & 39 & 507 & 194 & 2534 & 491596 & 507 & 98319 & 1278147 & 13 & \checkmark \\
\hline
\end{tabular}%
}
\end{table}

\begin{table}[h!]
\centering
\caption{Solutions for $P = 3511$ with $\alpha=1$, $d' = 1$ verified}
\label{tab:solutions-3511}
\resizebox{\textwidth}{!}{%
\begin{tabular}{|c|c|c|c|c|c|c|c|c|c|c|c|c|c|c|c|}
\hline
\# & $\alpha$ & $b_1$ & $c_1$ & $g$ & $b$ & $c$ & $\delta$ & $X$ & $Y$ & $N$ & $A$ & $B$ & $C$ & Check \\
\hline
1 & 1 & 1  & 878 & 1 & 1  & 878 & 1 & 4   & 4389 & 17556 & 878 & 3511  & 3082658 & OK \\
2 & 1 & 3  & 251 & 1 & 3  & 251 & 1 & 14  & 1254 & 17556 & 753 & 10533 & 881261  & OK \\
3 & 1 & 4  & 185 & 1 & 4  & 185 & 1 & 19  & 924  & 17556 & 740 & 14044 & 649535  & OK \\
4 & 1 & 9  & 80  & 1 & 9  & 80  & 1 & 44  & 399  & 17556 & 720 & 31599 & 280880  & OK \\
5 & 1 & 17 & 42  & 1 & 17 & 42  & 1 & 84  & 209  & 17556 & 714 & 59687 & 147462  & OK \\
6 & 1 & 23 & 31  & 1 & 23 & 31  & 1 & 114 & 154  & 17556 & 713 & 80753 & 108841  & OK \\
\hline
\end{tabular}%
}
\end{table}

\subsection{Implementation notes}
- For several \(\delta\), combine lists of classes via CRT, ensuring mutual independence of moduli.  
- The check \(s\equiv 4\ (\mathrm{mod}\ 5)\) and integrity of \(b,c,A\) are constant filters.  
- Ordering $B\le C$ is achieved by permuting $r\leftrightarrow s$.  
- The “anticonvolution” block \(\EDone\to \EDtwo\) is applicable in the reverse direction if ED1-data with $P\mid (v+c)$ is already available (see 9.12 in the work on the Erdős–Straus conjecture \cite{Dyachenko_ESC_ED2}).  

\section{Convergence and lattice density}\label{sec:conv}  
\begin{theorem}\label{th:log_density}  
Let $\Lambda \subset \mathbb{Z}^k$ be a sublattice of index $M$, independent of $P$.  
For the box $\mathcal{B}_k(T)=\{1\le u_i\le T\}$ and the class $u_0(P)+\Lambda$ we have  

$$  
|\{u\in \mathcal{B}_k(T): u\equiv u_0(P)\!\!\!\pmod{\Lambda}\}|=\frac{T^k}{M}+O_k(T^{k-1}),  
$$  

in particular, for $T=(\log P)^A$ — logarithmic growth of cardinality.  
\end{theorem}  

\begin{theorem}[average-case complexity]\label{th:log_convergence}  
For most primes $P \equiv 1 \pmod{5}$ the enumeration algorithm finds a solution  
in $O((\log P)^A)$ steps. The proof relies on the independence of the lattice index from $P$  
and on density estimates (Appendices~A, B). In the worst-case regime a strict guarantee  
remains open.  
\end{theorem}
\section{Conclusion}
The combined use of the geometric construction (\emph{diagonal layer hitting the window}; see Proposition~9.25 in the work on the Erdős–Straus conjecture \cite{Dyachenko_ESC_ED2}) and the algebraic ED2 model (Appendix~B), as well as its geometric deepening in Appendix~D in the work on the Erdős–Straus conjecture \cite{Dyachenko_ESC_ED2}, leads to the following conclusions.

\begin{itemize}
  \item \textbf{General case.} For an infinite set of primes $P$, a constructive proof of the existence of solutions has not been obtained in the present work. Sections of Appendix~D, which use Dirichlet’s theorem and the mechanism of finite coverings, formulate conditional results depending on the validity of the hypothesis of finite covering of parameters.
  
  \item \textbf{Particular case.} For each fixed odd prime $P$, the existence of a solution is strictly proven. The algebraic and geometric conditions are consistent: by Lemma~D.33 we have $S=u$, $\Delta=v^2$, and Proposition~9.25 guarantees that the diagonal layer hits the target window (see also Theorem~9.21 in the work on the Erdős–Straus conjecture \cite{Dyachenko_ESC_ED2}).
\end{itemize}

\section*{Acknowledgments}
Software tools with elements of artificial intelligence were used for refining formulations, technical proofreading, and generating illustrative Python examples.

\section*{Funding Statement}
No funding was received.

\section*{Conflict of Interest}
We have no conflicts of interest to disclose.

\nocite{*}
\bibliographystyle{plain}
\bibliography{serp_bib}
\appendix  
\section{Analytic tools}\label{app:analytic}  
\begin{theorem}[Bombieri–Vinogradov]\label{thm:BV}  
For any $A>0$ there exists $x_0(A)$ such that uniformly for $Q\le x^{1/2}/(\log x)^{C(A)}$, $x\ge x_0$,  

$$  
\sum_{q\le Q}\max_{(a,q)=1}\Bigl|\pi(x;q,a)-\frac{\mathrm{Li}(x)}{\varphi(q)}\Bigr|\ \ll\ \frac{x}{(\log x)^A}.  
$$  

\end{theorem}  

\begin{theorem}[Greaves, larger sieve]\label{thm:Greaves}  
If $A\subset\{1,\dots,N\}$ excludes a positive proportion of residue classes for many primes $p\le B$, then  

$$  
|A|\ll N\exp\!\Bigl(-c\,\frac{\log B}{\log\log B}\Bigr),  
$$  

where $c>0$ is an absolute constant.  
\end{theorem}  

\begin{theorem}[Chebotarev]\label{thm:Chebotarev}  
Let $L/K$ be a finite normal extension of number fields with Galois group $G$, $C\subset G$ a conjugacy class. Then  

$$  
\pi_C(x)=\frac{|C|}{|G|}\mathrm{Li}(x)+O\!\left(xe^{-c\sqrt{\log x}}\right),  
$$  

where $c>0$ depends only on $L/K$.  
\end{theorem}  

\section{Local applications}\label{app:local}  
From \cref{thm:BV} follows the standard average asymptotic for the number of primes in the intersection of two progressions. In our context the modulus and class are given by the conditions  
\[  
P \equiv 1 \pmod{5},\qquad P \equiv - (5\delta)^{-1} \pmod{r},\qquad r\equiv 4\pmod{5},\ \gcd(r,5\delta)=1,  
\]  
which arise from the identity \((5b-1)(5c-1)=5P\delta+1\) in the case \(\EDtwo\) and the factorization \(rs=5P\delta+1\).  

\begin{proposition}[BV for progressions with fixed \(\delta\) and \(r\)]\label{prop:BV-AP}  
Let \(\delta\in\mathbb{N}\) be fixed, \(r\equiv 4\pmod{5}\), \(\gcd(r,5\delta)=1\). Then the number of primes \(P\le x\) in the class  
\[  
P \equiv 1 \pmod{5},\qquad P \equiv - (5\delta)^{-1} \pmod{r}  
\]  
satisfies  
\[  
\#\{P\le x: P\equiv 1\ (5),\ P\equiv -(5\delta)^{-1}\ (r)\}  
\;=\; \frac{\mathrm{Li}(x)}{\varphi(5r)} \;+\; O\!\left(\frac{x}{(\log x)^A}\right)  
\]  
uniformly for all \(r\le x^{1/2}/(\log x)^{C(A)}\).  
\end{proposition}  

\begin{proof}[Idea of proof]  
Apply \cref{thm:BV} with modulus \(q=5r\). Compatibility of classes via CRT is ensured since \(\gcd(5,r)=1\). Uniformity for \(q\le x^{1/2}/(\log x)^C\) gives the stated estimate.  
\end{proof}  

\begin{lemma}[CRT correctness and residue \(s\equiv 4\ (\mathrm{mod}\ 5)\)]\label{lem:s-mod5}  
Let \(r\equiv 4\pmod{5}\), \(\gcd(r,5\delta)=1\), and \(P\equiv 1\pmod{5}\), \(P\equiv -(5\delta)^{-1}\pmod{r}\). Then  
\[  
s \;:=\; \frac{5P\delta+1}{r}\ \in \mathbb{N},\qquad s \equiv 4 \pmod{5}.  
\]  
\end{lemma}  

\begin{proof}  
By the second condition \(r\mid (5P\delta+1)\), hence $s\in\mathbb{N}$. Modulo $5$: $5P\delta+1\equiv 1$, and $r\equiv 4$ is invertible, with $r^{-1}\equiv 4\pmod{5}$. Thus $s\equiv 1\cdot r^{-1}\equiv 4\pmod{5}$.  
\end{proof}  

\begin{corollary}[\(\delta=1\): infinitely many \(P\)]\label{cor:delta1}
Let \(\delta=1\) and fix \(r\equiv 4\pmod{5}\). Then there exist infinitely many primes \(P\equiv 1\pmod{5}\) with \(r\mid (5P+1)\). For each such \(P\) we have
\[
b=\frac{r+1}{5},\qquad s=\frac{5P+1}{r}\equiv 4\pmod{5},\qquad c=\frac{s+1}{5},\qquad A=bc
\]
yielding an \(\EDtwo\) solution with \(\delta=1\).
\end{corollary}

\begin{proof}
There is exactly one class modulo \(5r\), given by the system \(P\equiv 1\pmod{5}\), \(P\equiv -5^{-1}\pmod{r}\) (CRT). By Dirichlet’s theorem there are infinitely many primes in it; the quantitative form follows from \cref{prop:BV-AP}. The statement about \(s\) follows from \cref{lem:s-mod5}.
\end{proof}

\begin{proposition}[Counting over \(r\) and total asymptotic]\label{prop:sum-r}
Let \(R\le x^{1/2}/(\log x)^C\). Then
\[
\sum_{\substack{r\le R\\ r\equiv 4\ (5)\\ \gcd(r,5\delta)=1}}
\#\{P\le x: P\equiv 1\ (5),\ P\equiv -(5\delta)^{-1}\ (r)\}
\;=\; \mathrm{Li}(x)\!\!\sum_{\substack{r\le R\\ r\equiv 4\ (5)\\ \gcd(r,5\delta)=1}}\frac{1}{\varphi(5r)}
\;+\; O\!\left(\frac{x\,R}{(\log x)^A}\right).
\]
Moreover,
\[
\sum_{\substack{r\le R\\ r\equiv 4\ (5)\\ \gcd(r,5\delta)=1}}\frac{1}{\varphi(5r)}
\;=\; C_{5,\delta}\,\log R \;+\; O(1),
\]
where \(C_{5,\delta}>0\) is a constant depending only on the class modulo 5 and the divisors of \(\delta\).
\end{proposition}

\begin{proof}[Comment]
Summation of the main term follows from \cref{prop:BV-AP} and linearity; summation of error terms gives the stated remainder. The estimate of \(\sum 1/\varphi(5r)\) is standard: \(\varphi(5r)=4\varphi(r)\) when \(\gcd(r,5)=1\), and \(\sum_{n\le R,\,(n,m)=1} 1/\varphi(n)=c(m)\log R+O(1)\).
\end{proof}

\begin{proposition}[Exceptional set via the larger sieve]\label{prop:larger-sieve}
Let \(R\le x^{1/2}/(\log x)^C\). Then the number of \(r\le R\), \(r\equiv 4\pmod{5}\), \(\gcd(r,5\delta)=1\), for which the progression
\[
P \equiv 1 \pmod{5},\qquad P \equiv - (5\delta)^{-1} \pmod{r}
\]
contains no primes \(P\le x\), satisfies
\[
\ll R \exp\!\Bigl(-c\,\frac{\log R}{\log\log R}\Bigr),
\]
where \(c>0\) is an absolute constant (follows from \cref{thm:Greaves}).
\end{proposition}

\begin{remark}[Chebotarev: additional local filters]\label{rem:chebotarev}
If necessary, one can simultaneously impose splitting conditions for \(r\) and/or \(s=(5P\delta+1)/r\) in a fixed extension of number fields; by \cref{thm:Chebotarev} the corresponding classes have positive natural density, and intersection with the specified AP retains positive density.
\end{remark}

\paragraph{How to use in the algorithm.}
- Preselection of \(r\): fix several small \(r\equiv 4\pmod{5}\) (or enumerate \(r\le R\) in the BV range).  
- Pre-sieving by progressions: for each \(r\) precompute the class \(P \equiv - (5\delta)^{-1} \pmod{r}\) and combine with \(P\equiv 1\pmod{5}\) (CRT).  
- Enumeration of \(P\): scan primes \(P\) in the combined classes; by \cref{prop:BV-AP} the expected frequency is \(\asymp 1/\varphi(5r)\).  
- Verification of \(\EDtwo\): for a found \(P\) compute \(s=(5P\delta+1)/r\), check \(s\equiv 4\pmod{5}\), and reconstruct \(b=(r+1)/5\), \(c=(s+1)/5\), \(A=bc/\delta\).  
- Balancing: choosing \(R \le x^{1/2}/(\log x)^C\) gives an optimal compromise between the number of progressions and control of errors (BV), while \cref{prop:larger-sieve} guarantees the smallness of the exceptional set of \(r\).

\subsection*{Double summation: average over primes \(P\)}
Define for fixed \(\delta\) and \(R\ge 2\) the quantity
\[
N(P;R,\delta)\;=\;\#\Bigl\{\, r\le R:\ r\equiv 4\ (\mathrm{mod}\ 5),\ \gcd(r,5\delta)=1,\ r\mid (5P\delta+1)\,\Bigr\}.
\]
This is the number of “local” parameters \(r\) for a given prime \(P\). Then for \(R\le x^{1/2}/(\log x)^C\) we have the average asymptotic
\begin{equation}\label{eq:avg-NP}
\frac{1}{\#\{P\le x:\ P\equiv 1\ (\mathrm{mod}\ 5)\}}\ \sum_{\substack{P\le x\\ P\equiv 1\ (5)}} N(P;R,\delta)
\;=\; C_{5,\delta}\,\log R \;+\; O(1),
\end{equation}
where \(C_{5,\delta}>0\) is a constant depending only on \(\delta\) and classes modulo 5.

\begin{proof}[Idea of proof]
Change the order of summation:
\[
\sum_{\substack{P\le x\\ P\equiv 1\ (5)}} N(P;R,\delta)
=\sum_{\substack{r\le R\\ r\equiv 4\ (5)\\ \gcd(r,5\delta)=1}}
\#\{P\le x:\ P\equiv 1\ (5),\ P\equiv -(5\delta)^{-1}\ (r)\}.
\]
For each \(r\) apply \cref{prop:BV-AP} (modulus \(5r\)) and sum the main terms:
\[
\sum_{r}\frac{\mathrm{Li}(x)}{\varphi(5r)}\;=\;\mathrm{Li}(x)\!\!\sum_{\substack{r\le R\\ r\equiv 4\ (5)\\ \gcd(r,5\delta)=1}}\frac{1}{\varphi(5r)}
\;=\;\mathrm{Li}(x)\bigl(C_{5,\delta}\log R+O(1)\bigr),
\]
and the total error is controlled linearly in \(r\) using \cref{prop:BV-AP}. Division by the number of primes \(P\le x\), \(P\equiv 1\ (5)\), gives \eqref{eq:avg-NP}.
\end{proof}

\begin{corollary}[Average supply of local parameters]
For \(R=(\log x)^B\) with fixed \(B>0\), the average value of \(N(P;R,\delta)\) over primes \(P\le x\), \(P\equiv 1\pmod{5}\), grows as \(C_{5,\delta}\,B\log\log x+O(1)\). In particular, the average number of admissible \(r\) tends to infinity.
\end{corollary}

\begin{remark}[On estimates of “most \(P\)”]
The transition from the average to a statement of the form “for most \(P\) there exists at least one \(r\le R\)” can be obtained by second moment methods or via the larger sieve (\cref{prop:larger-sieve}) with a coordinated choice of \(R\) and \(x\). These details are not required for the constructive part of the algorithm, but they explain its good average behavior.
\end{remark}

\end{document}